\documentclass[12pt,a4paper]{article}
\usepackage{color}
\usepackage{amsmath, amssymb, theorem, latexsym}
\usepackage{graphicx}
\usepackage{epsf}
\usepackage{fullpage}

\newtheorem{theorem}{Theorem}[section]
\newtheorem{lemma}[theorem]{Lemma}
\newtheorem{proposition}[theorem]{Proposition}

\newtheorem{corollary}[theorem]{Corollary}

\allowdisplaybreaks
\newcommand{\beq}{\begin{equation}}
\newcommand{\eeq}{\end{equation}}

\def \cL {\mathcal L}

\numberwithin{equation}{section}

\newcommand{\noib}{\noindent $\bullet $~}
\newcommand{\ie}{\emph{i.e. }}

\newcommand{\resp}{\emph{resp. }}

\newcommand{\piq}{\frac{\pi}{4}}

\newcommand{\hH}{\hat{H}}

\def \R {{\mathbb R}}

\newcommand{\cE}{\mathcal{E}}

\title{On the number of nodal domains of the
2D isotropic quantum harmonic oscillator\\
-- an extension of results
of A.~Stern--}

\author{P. B\'erard \\ Institut Fourier, Universit\'{e} de Grenoble and CNRS, B.P.74,
\\ F 38402 Saint Martin d'H\`{e}res Cedex, France.\\
and \\
B. Helffer \\
Laboratoire de Math\'ematiques, Univ. Paris-Sud 11 and CNRS,\\
F 91405 Orsay Cedex, France, and\\
Laboratoire Jean Leray, Universit\'{e} de Nantes.}

\date{September 8, 2014}

\begin{document}
\maketitle

\begin{abstract} In the case of the sphere and the square, Antonie
Stern (1925)  claimed  in her  PhD thesis the existence of an
infinite sequence of eigenvalues whose corresponding eigenspaces
contain an eigenfunction with two nodal domains. These two
statements were  given complete proofs respectively by Hans Lewy in
1977, and the authors in 2014 (see also
Gauthier-Shalom--Przybytkowski (2006)). The aim of this paper is to
obtain a similar result in the case of the isotropic  quantum
harmonic oscillator in the two dimensional case.
\end{abstract}

Keywords: Quantum harmonic oscillator, Nodal lines, Nodal domains,
Courant theorem.\\
MSC 2010: 35B05, 35Q40, 35P99, 58J50, 81Q05.

\section{Introduction and main results}

The aim  of this paper is to construct a sequence of
eigenvalues, and a corresponding sequence of eigenfunctions,
for the  2D isotropic quantum harmonic oscillator
\begin{equation}
\widehat
H:= -\Delta + x^2 + y^2\,,
\end{equation} with exactly two nodal domains. Similar
results were stated  in Antonie Stern' PhD thesis (with R. Courant
as advisor), for the Dirichlet problem for the Laplacian in the
square, and in the case of the Laplace-Beltrami  operator on the
sphere $\mathbb S^2$ \cite{St}. This is only much later that H. Lewy
\cite{Lew} in the case of the sphere (see also \cite{BeHe2}), and
the authors \cite{BeHe} in the case of the square with Dirichlet
conditions (see also Gauthier-Shalom--Przybytkowski \cite{GSP}),
have proposed complete proofs of her statements.

Coming back to the isotropic harmonic oscillator,   an
orthogonal basis of eigenfunctions is given by
\begin{equation}
\phi_{m,n} (x,y) = H_m(x) H_n (y) \exp( - \frac{x^2+y^2}{2})\,,
\end{equation}
 for $(m,n)\in \mathbb N^2$, where $H_n(x)$ denotes the
Hermite polynomial of degree $n$.

The eigenfunction $\phi_{m,n}$ corresponds to the eigenvalue
$2(m +n+1)\,,$
\begin{equation}
\widehat H \phi_{m,n} = 2(m+n+1)\, \phi_{m,n}\, .
\end{equation}
Here we use the definitions and notation of Szeg\"o \cite[\S
5.5]{Sz}.

The eigenspace $\cE_{\ell}$ associated with the eigenvalue
$\hat{\lambda}(\ell)= 2(\ell +1)$ has dimension $(\ell + 1)$, and is
generated by the eigenfunctions $\phi_{\ell,0}$, $\phi_{\ell-1,1}$,
\dots, $\phi_{0,\ell}$.

For $\theta \in [0,\pi[$, we shall consider the families of
eigenfunctions,
\begin{equation}
\Phi^{\theta}_{n} := \cos\theta \, \phi_{n,0} + \sin\theta \,
\phi_{0,n}\,,
\end{equation}
corresponding to the eigenvalue $2(n+1)$.

Our aim is to prove the following  theorems.

\begin{theorem}\label{prop1}   Assume that  $n$ is odd. Then,
there exists an open interval $I_{\frac{\pi}{4}}$ containing
$\frac{\pi}{4}$, and an open interval $I_{\frac{3\pi}{4}}$,
containing $\frac{3\pi}{4}$, such that for
$$
\theta \in I_{\frac{\pi}{4}} \cup I_{\frac{3\pi}{4}} \setminus
\{\frac{\pi}{4},\frac{3\pi}{4}\}\,,
$$
the nodal set $N(\Phi^{\theta}_n)$ is a connected simple regular
curve, and the eigenfunction $\Phi^{\theta}_n$ has two nodal domains
in $\R^2$.
\end{theorem}

\begin{theorem}\label{th2}
 Assume that $n$ is odd. Then, there exists $\theta_c >0$ such
that, for $0 < \theta < \theta_c$, the nodal set
$N(\Phi^{\theta}_n)$ is a connected simple regular curve, and the
eigenfunction $\Phi^{\theta}_n$ has two nodal domains in $\R^2$.
\end{theorem}

As in the case of the square,  to prove Theorem~\ref{prop1}, we
begin by a symmetry argument to reduce to a neighborhood of either
$\pi/4$ or $3\pi/4$, say $3\pi/4$. The first step in the proof is to
analyze the zero set when $\theta = \frac{3\pi}{4}$, in particular
the points where the zero set is critical, and to show that this
only occurs on the diagonal.

The second  step is then to show that the double crossings on
the diagonal disappear by perturbation, for $\theta$ close to
and different from $\frac{3\pi}{4}$\,. Using the local nodal
patterns and some barrier lemmas, one can then show that the nodal
set becomes a connected simple curve, asymptotic to $x=y$ at $\pm
\infty\,$. The local stability of the nodal set under perturbation,
then gives an explicit interval containing $3 \pi/4$ in which the
phenomenon occurs (for $\theta \neq \frac{3\pi}{4}$ of course).

The proof of Theorem~\ref{th2} follows similar lines.


\newpage

\textbf{Remarks}.\vspace{-3mm}
\begin{enumerate}
\item In the case of the square the same kind of analysis is also
interesting for determining when the number of nodal domains of an
$n$-th eigenfunction  is equal to $n$ (the so called Courant sharp
situation). Pleijel \cite{Pl} observed that in the case of the
square (Dirichlet condition) this  only occurs for the first, the
second and the fourth eigenfunctions (see \cite{BeHe} for a complete
argument). In the case of the sphere, as a consequence of the
analysis of Leydold \cite{Ley}, this  only occurs for the first and
the second eigenfunctions. It is natural to investigate a similar
question in the case of the isotropic harmonic oscillator, and more
generally, the validity of Pleijel's theorem in this case.  In
the last section, we will give a Leydold's like proof of the fact
that the only Courant sharp eigenvalues of the harmonic oscillator
are $\hat{\lambda }(\ell)$, for $\ell = 0, 1$ and $2$. As
communicated by I. Polterovich, this question will be analyzed from
a different point of view in \cite{Cha}.
\item A connected question is to analyze the zero set when $\theta$
is a random variable. We refer to \cite{HZZ} for  results in
this direction.
\item These questions are related to the question of spectral minimal
partitions \cite{HHOT}. In the case of the harmonic oscillator
similar questions appear in the analysis of the properties of
ultracold atoms (see for example \cite{RL}).
\end{enumerate}

 Theorems~\ref{prop1} and \ref{th2} concern the eigenspace
$\cE_n$ of the harmonic oscillator, with $n$ odd. When $n$ is even,
the picture is different. Some nodal sets have compact connected
components, with or without critical zeros, some have both compact
and non-compact components. Other examples can be analyzed as well.
This will be analyzed in the future.

{\bf Acknowledgements.}\\
The second author would like to thank D. Jakobson, I. Polterovich
and M. Persson-Sundqvist for useful discussions, transmission of
information or computations.

\section{A  reminder on Hermite polynomials}

We use the definition, normalization, and notations of Szeg\"{o}'s book
\cite{Sz}. With these choices, $H_n$ has the following properties,
\cite[\S~5.5 and Theorem~6.32]{Sz}.\vspace{-3mm}

\begin{enumerate}
\item  $H_n$ satisfies the differential equation
$$
y''(t) - 2 t \,y'(t) + 2n\, y(t) = 0\,.
$$
\item  $H_n(t)$ is a polynomial of degree $n$ which is even (\resp odd) for $n$
even (\resp odd).
\item  $H_n(t) = 2 t\, H_{n-1}(t) - 2 (n-1) \,
H_{n-2}(t)\,, ~n\ge 2\,, ~~H_0(t)= 1\,, ~~H_1(t) = 2t\,.$
\item $H_n$ has $n$ simple zeros $t_{n,1} < t_{n,2} < \cdots <
t_{n,n}\,$.
\item
$$
 H_n(t) =  2 t\, H_{n-1}(t) - H'_{n-1}(t)\,.
$$
\item \begin{equation}
H'_n(t) = 2n H_{n-1} (t)\,.
\end{equation}
\item The coefficient of $t^n$ in  $H_n$ is $2^n$.
\item
$$
\int _{-\infty}^{+\infty} e^{-t^2} |H_n(t)| ^2\, dt = \pi^\frac 12 \,2^n\, n! \,.
$$
\item The  first zero $t_{n,1}$ of $H_n$ satisfies
\begin{equation}
 t_{n,1} = (2n+1)^\frac 12 - 6^{-\frac 12} (2n+1)^{-\frac 16}
(i_1+\epsilon_n)\,,
\end{equation}
where $i_1$ is the first positive real zero of the Airy function,
and $\lim_{n\rightarrow +\infty} \epsilon_n =0\,$.
\end{enumerate}

The following result (Theorem 7.6.1 in Szeg\"o's book \cite{Sz}) will
also be useful:

\begin{lemma}\label{L-leh}
The successive relative maxima of $|H_n(t)|$ form an increasing
sequence for $t \geq 0$\,.
\end{lemma}
{\bf Proof.}
\\
It is enough to observe that the function
$$
\Theta_n (t):= 2 n H_n(t)^2 + H_n'(t)^2
$$
satisfies
$$
\Theta_n'(t) = 4 t \, (H'_n(t)) ^2 \,.
$$\hfill $\square$

\section{Stern-like constructions for the harmonic oscillator, case $n$-odd}

\subsection{The case of the square}

 Consider the square  $[0,\pi]^2$, with Dirichlet boundary
conditions, and the following families of eigenfunctions associated
with the eigenvalues $\hat{\lambda}(1,2r) := 1+ 4r^2$, where $r$ is
a positive integer, and $\theta \in [0,\pi/4]$,
$$
(x,y) \mapsto \cos \theta \, \sin x \, \sin(2ry)  + \sin \theta \,
\sin(2rx) \, \sin y\,.
$$

\begin{figure}[!ht]
\begin{center}
\includegraphics[width=14cm]{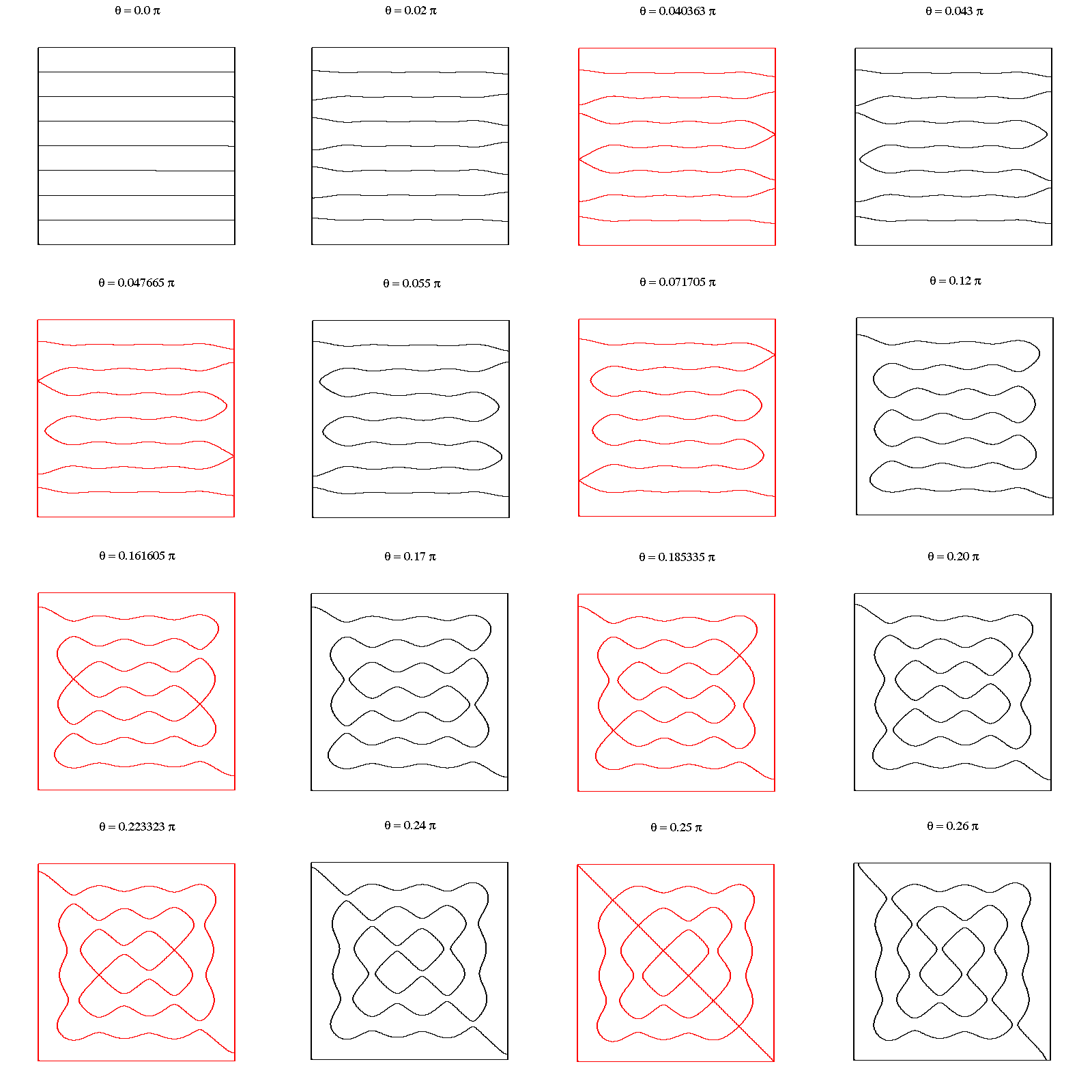}
\caption{Evolution of the nodal set in the case of the square.} \label{sq-1}
\end{center}
\end{figure}

According to \cite{St}, for any given $r$, the typical evolution of
the nodal sets when $\theta$ varies is similar to the case $r=4$
shown in Figure~\ref{sq-1} \cite[Figure~6.9]{BeHe}: generally
speaking, the nodal sets deform continuously, except for finitely
many values of $\theta$, for which crossings appear or disappear.

We would like to get similar results for the isotropic quantum
harmonic oscillator.

\subsection{Symmetries}

 Recall the notation,
\begin{equation}\label{P2}
\Phi^{\theta}_n(x,y) :=: \Phi_n(x,y,\theta) := \cos\theta \,
\phi_{n,0} + \sin\theta \, \phi_{0,n}\,.
\end{equation}
Since $\Phi^{\theta + \pi}_n = - \Phi^{\theta}_n$, it suffices to
vary the parameter $\theta$ in the interval $[0,\pi[$.

\emph{Since $n$ is odd}, we have the following symmetries.

\begin{equation}\label{P4}
\left\lbrace
\begin{array}{ll}
\Phi^{\theta}_n(-x,y) & = \Phi^{\pi-\theta}_n(x,y)\,,\\[4pt]
\Phi^{\theta}_n(x,-y) & = - \Phi^{\pi-\theta}_n(x,y)\,,\\[4pt]
\Phi^{\theta}_n(y,x) & = \Phi^{\frac{\pi}{2}-\theta}_n(x,y)\,.\\[4pt]
\end{array}
\right.
\end{equation}

When $n$ is odd, it therefore suffices to vary the parameter
$\theta$ in the interval $[0,\piq]$. The case $\theta = 0$ is
particular, so that we shall mainly consider $\theta \in ]0,\piq]$.

\subsection{Critical zeros}

A \emph{critical zero} of $\Phi^{\theta}_n$ is a point $(x,y) \in
\R^2$ such that both $\Phi^{\theta}_n$ and its differential
$d\Phi^{\theta}_n$ vanish at $(x,y)$. The critical zeros of
$\Phi^{\theta}_n$ satisfy the following equations.

\begin{equation}\label{P6}
\left\lbrace
\begin{array}{ll}
\cos\theta \, H_n(x) + \sin\theta \, H_n(y) &= 0\,,\\[4pt]
\cos\theta \, H'_n(x) & = 0 \,,\\[4pt]
\sin\theta \, H'_n(y) &= 0\,.\\
\end{array}
\right.
\end{equation}

Equivalently, using the properties of the Hermite polynomials, a
point $(x,y)$ is a critical zero of $\Phi^{\theta}_n$ if and only if

\begin{equation}\label{P8}
\left\lbrace
\begin{array}{ll}
\cos\theta \, H_n(x) + \sin\theta \, H_n(y) &= 0\,,\\[4pt]
\cos\theta \, H_{n-1}(x) & = 0 \,,\\[4pt]
\sin\theta \, H_{n-1}(y) &= 0\,.\\
\end{array}
\right.
\end{equation}

The only possible critical zeros of the eigenfunction
$\Phi^{\theta}_n$ are the points $(t_{n-1,i}\,,t_{n-1,j})$ for $1
\le i,j \le (n-1)$, where the coordinates are the zeros of the
Hermite polynomial $H_{n-1}$. The point $(t_{n-1,i}\,,t_{n-1,j})$ is
a critical zero of $\Phi^{\theta}_n$ if and only if $\theta =
\theta(i,j)$, where $\theta(i,j) \in ]0,\pi[$ is uniquely defined by
the equation,
\begin{equation}\label{P10}
\cos\left( \theta(i,j) \right) \, H_n(t_{n-1,i}) + \sin\left(
\theta(i,j) \right) \, H_n(t_{n-1,j}) = 0\,.
\end{equation}
Here we have used the fact that $H_n$ and $H'_n$ have no common
zeros. We have proved the following lemma.

\begin{lemma}\label{L-cz}
For $\theta \in [0,\pi[$, the eigenfunction $\Phi^{\theta}_n$ has no
critical zero, unless $\theta$ is one of the $\theta(i,j)$ defined
by equation \eqref{P10}. In particular $\Phi^{\theta}_n$ has no
critical zero, except for finitely many values of the parameter
$\theta \in [0,\pi[$. Let $\theta_0 = \theta(i_0,j_0)$,
defined by some $(t_{n-1,i_0}\,, t_{n-1,j_0})$. The function
$\Phi^{\theta_0}_n$ has finitely many critical zeros, namely the
points $ (t_{n-1,i}\,, t_{n-1,j})$ which satisfy
\begin{equation}
\cos\theta_0 \, H_n(t_{n-1,i}) + \sin\theta_0 \, H_n(t_{n-1,j}) =
0\,,
\end{equation}
among them the point  $(t_{n-1,i_0}\,, t_{n-1,j_0})$.
\end{lemma}

\textbf{Remarks}.\\
 From the general properties of nodal lines
\cite[Properties~5.2]{BeHe}, we  derive the following facts.
\vspace{-3mm}
\begin{enumerate}
    \item When $\theta \not \in \left\lbrace \theta(i,j) ~|~ 1 \le i,j \le n-1
    \right\rbrace$, the nodal set of the eigenfunction $\Phi^{\theta}_n$, denoted by
    $N(\Phi^{\theta}_n)$, is a smooth $1$-dimensional submanifold of $\R^2$.
    \item When $\ \theta  \in \left\lbrace \theta(i,j)~|~ 1 \le i,j \le n-1
    \right\rbrace$, the nodal set $N(\Phi^{\theta}_n)$ has finitely
    many singularities which are double crossings. Indeed, the Hessian of the
    function $\Phi^{\theta}_n$ at a critical zero $(t_{n-1,i},t_{n-1,j})$ is
    given by
    \begin{equation*}\label{hessian}
    \mathrm{Hess}_{(t_{n-1,i},t_{n-1,j})}\Phi^{\theta}_n =
    \exp{(-\frac{t_{n-1,i}^2 + t_{n-1,j}^2}{2})}\,
    \begin{pmatrix}
      \cos\theta \, H''_n(t_{n-1,i}) & 0 \\
      0 & \sin\theta \, H''_n(t_{n-1,j}) \\
    \end{pmatrix}\,,
    \end{equation*}
and the assertion follows from the fact that $H_{n-1}$ has simple
zeros.
\end{enumerate}

\subsection{General properties of the nodal set $N(\Phi^{\theta}_n)$}

Denote by $\cL$ the finite lattice
\begin{equation}\label{P22}
\cL := \left\lbrace (t_{n,i}\,,t_{n,j}) ~|~ 1 \le i,j \le n
\right\rbrace \subset \R^2 \,,
\end{equation}
consisting of points whose coordinates are the zeros of the Hermite
polynomial $H_n$. Since we can assume that $\theta \in ]0,\piq]$, we
have the following inclusions for the nodal set,

\begin{equation}\label{P24}
\cL \subset N(\Phi^{\theta}_n) \subset \cL \cup \left\lbrace (x,y)
\in \R^2 ~|~ H_n(x) \, H_n(y) < 0 \right\rbrace \,.
\end{equation}

\textbf{Remarks}. Assume that $\theta \in ]0,\piq]$.\\
(i) The nodal set $N(\Phi^{\theta}_n)$ cannot meet the vertical
lines $\{x=t_{n,i}\}$, or the horizontal lines $\{y=t_{n,i}\}$ away
from the set $\cL$. \\
(ii) The lattice point $(t_{n,i},t_{n,j})$ is not a
critical zero of $\Phi^{\theta}_n$ (because $H_n$ and $H'_n$
have no common zero). As a matter of fact, near a lattice point, the
nodal set $N(\Phi^{\theta}_n)$ is a single arc through the lattice
point, with a tangent which is neither horizontal, nor
vertical.\medskip

\begin{figure}[!ht]
\begin{center}
\includegraphics[width=10cm]{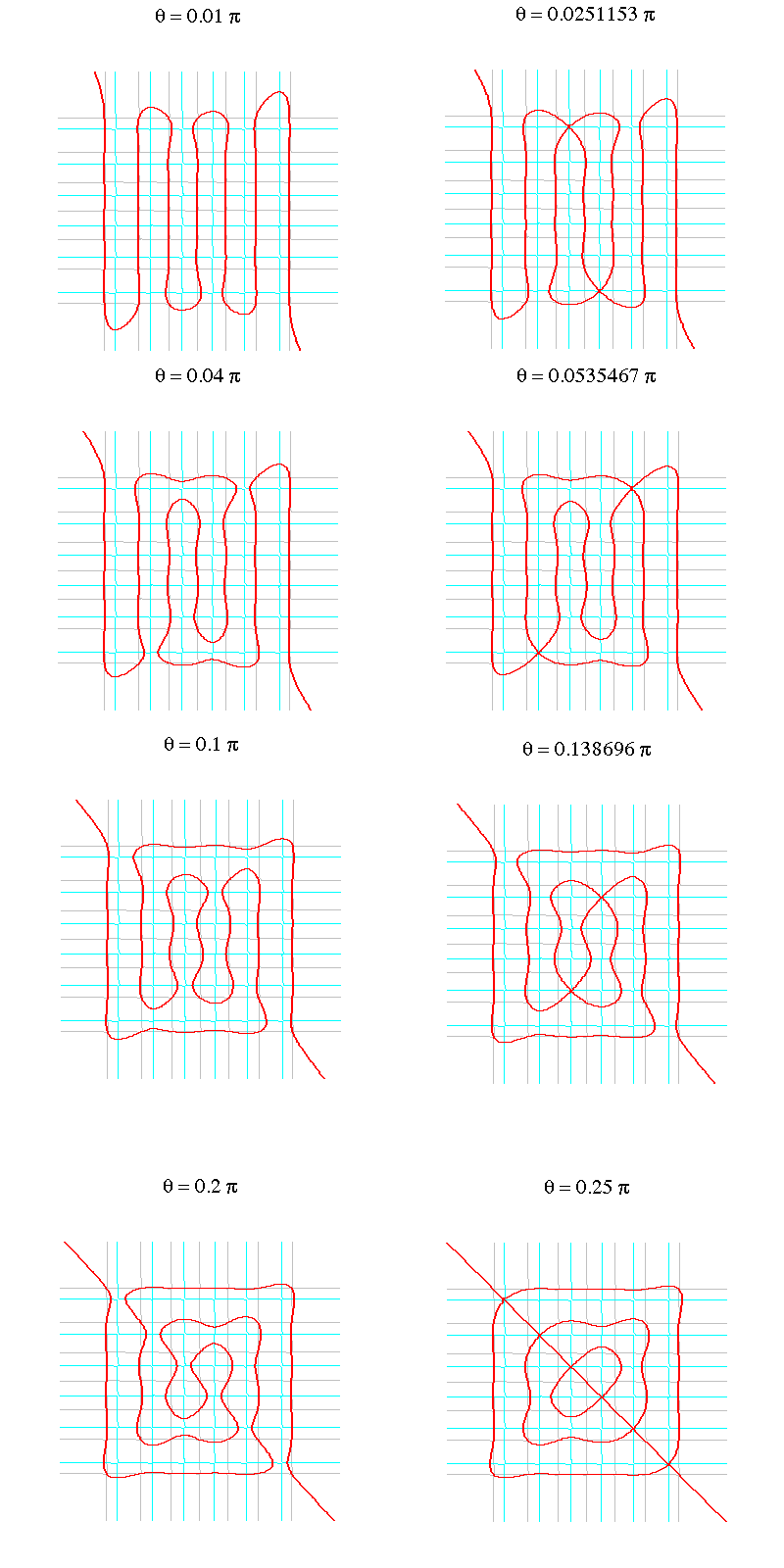}
\caption{Evolution of the nodal set $N(\Phi^{\theta}_n)$, for
$\theta \in ]0,\piq]$.} \label{FSch-1}
\end{center}
\end{figure}

Figure~\ref{FSch-1} shows the evolution of the nodal set of
$\Phi^{\theta}_n$ when $\theta$ varies in the interval $]0,\piq]$.
The values of $\theta$ with two digits are regular values (\ie
correspond to an eigenfunction without critical zeros), the values
of $\theta$ with at least three digits are critical values (\ie
correspond to an eigenfunction with critical zeros). The form of the
nodal set is stable between two consecutive critical values of the
parameter $\theta$. In the figures, the grey lines correspond to the
zeros of $H_7$. The blue lines correspond to the zeros of $H'_7$,
\ie to the zeros of $H_6$.


\medskip

We now describe the nodal set $N(\Phi^{\theta}_n)$ outside a large
enough square which contains the lattice $\cL$.  For this
purpose, we give two \emph{barrier lemmas}.

\begin{lemma}\label{L-ext1}
Assume that $\theta \in ]0,\piq]$.  For $n\geq 1$, define
$t_{n-1,0}$ to be the unique point in $]-\infty, t_{n,1}[$ such that
$H_n(t_{n-1,0}) = - H_n(t_{n-1,1})$. Then,\vspace{-3mm}
\begin{enumerate}
    \item $\forall t \le t_{n,1}$, the function $y \mapsto
    \Phi^{\theta}_n(t,y)$ has exactly one zero in the interval
    $[t_{n,n},+\infty[$\,;
    \item $\forall t < t_{n-1,0}$, the function $y \mapsto
    \Phi^{\theta}_n(t,y)$ has exactly one zero in the interval\break
    $]-\infty,+\infty[$\,.
\end{enumerate}
\end{lemma}

\textbf{Proof}.\\
Let $v(y) := \exp(\frac{t^2+y^2}{2})\, \Phi^{\theta}_n(t,y)$. In
$]t_{n,n},+\infty[$, $v'(y)$ is positive, and $v(t_{n,n})\le 0\,$.
The first assertion follows. The local extrema of $v$ occur at the
points $t_{n-1,j}\,$, for $1 \le j \le (n-1)\,$. The second assertion
follows from the definition of $t_{n-1,0}\,$, and from the
inequalities,
\begin{equation*}
\begin{split}
\cos\theta \, H_n(t) + \sin\theta \, H_n(t_{n-1,j}) & \le
\frac{1}{\sqrt{2}}\Big( H_n(t) - |H_n|(t_{n-1,j})\Big)\\
& < -\, \frac{1}{\sqrt{2}} \Big( H_n(t_{n-1,1}) -
|H_n|(t_{n-1,j})\Big) \le 0\,,
\end{split}
\end{equation*}
where we have used Lemma~\ref{L-leh}.\hfill $\square$\medskip

\textbf{Remark}. Using the symmetry with respect to the vertical
line $\{x=0\}$, one has similar statements for $t \ge t_{n,n}$ and
for $t > - t_{n-1,0}\,$.

\begin{lemma}\label{L-ext2} Let $\theta \in ]0,\piq]$. Define
$t_{n-1,n}^{\theta} \in ]t_{n,n}\,,\,\infty[$ to be the unique point
such that $\tan\theta \, H_n(t_{n-1,n}^{\theta}) = H_n(t_{n-1,1})$.
Then,\vspace{-3mm}
\begin{enumerate}
    \item $\forall t \ge t_{n,n}$, the function $x \mapsto
    \Phi^{\theta}_n(x,t)$ has exactly one zero in the interval
    $]-\infty,t_{n,1}]$\,;
    \item $\forall t > t_{n-1,n}^{\theta}$, the function $x \mapsto
    \Phi^{\theta}_n(x,t)$ has exactly one zero in the interval
    $]-\infty,\infty[$\,.
    \item For $\theta_2 > \theta_1$, we have
    $t_{n-1,n}^{\theta_2} < t_{n-1,n}^{\theta_1}\,$.
\end{enumerate}
\end{lemma}

\textbf{Proof}. Let $h(x) := \exp(\frac{x^2+t^2}{2})\,
\Phi^{\theta}_n(x,t)$. In the interval $]-\infty,t_{n,1}]$, the
derivative $h'(x)$ is positive, $h(t_{n,1}) > 0$, and
$\lim_{x\rightarrow  -\infty} h(x) =-\infty$, since $n$ is odd. The
first assertion follows. The local extrema of $h$ are achieved at
the points $t_{n-1,j}$. Using Lemma~\ref{L-leh}, for $t \geq
t_{n-1,n}^{\theta}\,$, we have the inequalities,
\begin{equation*}
\begin{split}
H_n(t_{n-1,j}) + \tan\theta \, H_n(t) &\ge
\tan\theta \, H_n(t_{n-1,n}^{\theta}) - |H_n(t_{n-1,j})|\\
& = H_n(t_{n-1,1}) - |H_n(t_{n-1,j})| \ge 0\,.
\end{split}
\end{equation*}
\hfill $\square$\medskip

\textbf{Remark}. Using the symmetry with respect to the horizontal
line $\{y=0\}$, one has similar statements for $t \le t_{n,1}$ and
for $t < - t_{n-1,n}^{\theta}\,$.\medskip

As a consequence of the above lemmas, we have the following
description of the nodal set far enough from $(0,0)$.

\begin{proposition}\label{P-ext}
Let $\theta \in ]0,\piq]$. In the set $\R^2 \setminus
]-t_{n-1,n}^{\theta},t_{n-1,n}^{\theta}[\times ]t_{n-1,0},
|t_{n-1,0}|[$, the nodal set $N(\Phi^{\theta}_n)$ consists of two
regular arcs. The first arc is a graph $y(x)$ over the
interval\break $]-\infty,t_{n,1}]$, starting from the point
$(t_{n,1},t_{n,n})$ and escaping to infinity with,
$$ \lim_{x\to -\infty}\frac{y(x)}{x} = - \sqrt[n]{\cot\theta}\,.$$
The second arc is the image of the first one under the symmetry with
respect to $(0,0)$ in $\R^2$.
\end{proposition}

\subsection{Local nodal patterns}\label{SS-lnp}

As in the case of the square, we study the possible local nodal
patterns taking into account the fact that  the nodal set contains the
lattice points $\cL$, can only intersect the connected components of
the set $\{H_n(x) \, H_n(y) < 0\}$, and consists of a simple arc at
the lattice points. The following figure summarized the possible
nodal patterns in the interior of the square
\cite[Figure~6.4]{BeHe},

\begin{figure}[!ht]
\begin{center}
\includegraphics[width=10cm]{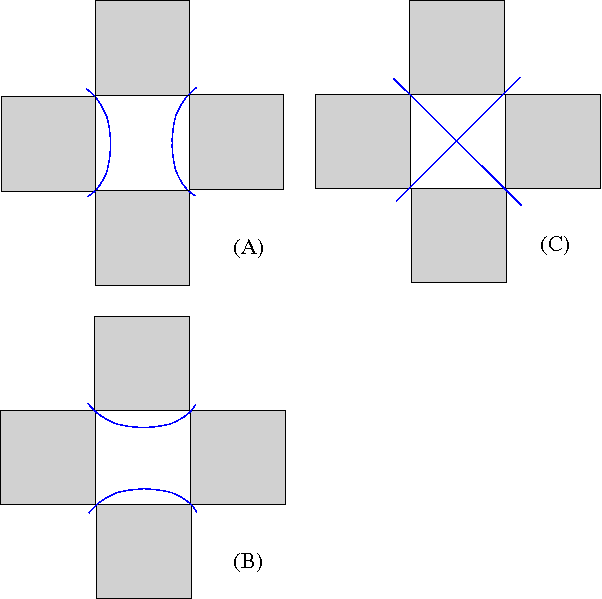}
\caption{Local nodal patterns.} \label{FSchr-local-nodal-patterns}
\end{center}
\end{figure}

Except for nodal arcs which escape to infinity, the local nodal
patterns for the quantum harmonic oscillator are the same (note that
in the present case, the connected components of the set $\{H_n(x)
\, H_n(y) < 0\}$ are rectangles, no longer equal squares).

Case (C) occurs near a critical zero. Following the same ideas as in
the case of the square, in order to decide between cases (A) and
(B), we use the barrier lemmas, Lemma~\ref{L-ext1} or~\ref{L-ext2},
the vertical lines $\{x=t_{n-1,j}\}$, or the horizontal lines
$\{y=t_{n-1,j}\}$.

\newpage

\section{Proof of Theorem \ref{prop1}}

Note that
$$\phi_{n,0}(x,y) - \phi_{0,n}(x,y) = - \Phi^{\frac{3\pi}{4}}_n(x,y) = -
\Phi^{\frac{\pi}{4}}_n(y,x)\,.
$$

Hence it is the same to work with $\theta =\frac \pi 4$ and the
anti-diagonal, or to work with $\theta = \frac{3\pi}{4}$ and the
diagonal. From now on, we work near $\frac{3\pi}{4}$.

\subsection{ The nodal set of $\Phi_n^{\frac{3\pi}{4}}$}

 \begin{proposition}\label{prop2} Let $\{t_{n-1,i}\,, 1 \le i \le
n-1\}$ denote the zeroes of $H_{n-1}$. For $n$ odd, the nodal set of
$\phi_{n,0} - \phi_{0,n}$ consists of the diagonal $x=y$, and of
$\frac{n-1}{2}$ disjoint simple closed curves crossing the diagonal
at the $(n-1)$ points $(t_{n-1,i}\,,t_{n-1,i})$, and the antidiagonal
at the $(n-1)$ points $(t_{n,i}\,, - t_{n,i})$.
\end{proposition}

To prove Proposition \ref{prop2}, we first observe that it is enough
to analyze the zero set of
$$
(x,y) \mapsto \Psi_{n}(x,y)  :=H_n(x) - H_n(y)\,.
$$
The critical points  of $\Psi_n$ are determined by
$$
H_n'(x)=0\,,\, H_n'(y)= 0\,.
$$

 Hence, the critical points of $\Psi_{n}$  consist of the
$(n-1)^2$ points $(t_{n-1,i}\,,t_{n-1,j})$, for \break $1 \le i,j \le (n-1)$,
where $t_{n-1,i}$ is the $i$-th zero of the polynomial $H_{n-1}\,$.

The zero set  of $\Psi_{n}$  contains the diagonal $\{x=y\}\,$.
 Since $n$ is odd, there  are only $n$ points belonging
to the zero set on the anti-diagonal $\{x+y=0\}$.

On the diagonal, there are $(n-1)$ critical points.  We claim
that there are no critical zeros outside the diagonal.
Indeed, let $(t_{n-1,i}\,,t_{n-1,j})$ be a critical zero. Then,
$H_n(t_{n-1,i}) = H_n(t_{n-1,j})$. Using Lemma~\ref{L-leh} and the
parity properties of Hermite polynomials, we see that
$|H_n(t_{n-1,i})| = |H_n(t_{n-1,j})|$ occurs if and only if
$t_{n-1,i} = \pm t_{n-1,j}\,$. Since $n$ is odd, we can conclude that
$H_n(t_{n-1,i}) = H_n(t_{n-1,j})$ occurs if and only if
$t_{n-1,i}=t_{n-1,j}\,$.

\subsection{ Existence of disjoint simple closed curves in the nodal
set of $\Phi_n^{\frac{3\pi}{4}}$}

The second part in the proof of the proposition follows closely the
proof in the case of the square (see Section 5 in \cite{BeHe}).
Essentially, the  Chebyshev polynomials are replaced by the Hermite
polynomials. Note however that the checkerboard is no more with
equal squares, and that  the square $[0,\pi]^2$ has to be
replaced in the argument by the rectangle  $[t_{n-1,0}, -t_{n-1,0}] \times
[-t_{n-1,n}^{\theta},t_{n-1,n}^{\theta}]$, for some $\theta$ such
that $0 < \theta < \frac{3\pi}{4}$, see Lemmas~\ref{L-ext1} and~\ref{L-ext2}.

The checkerboard argument holds, see \eqref{P24}.

The separation lemmas of our previous  paper \cite{BeHe} must
be substituted by Lemmas~\ref{L-ext1} and~\ref{L-ext2}, and similar
statements with the lines $\{x=t_{n-1,j}\}$ and $\{y=t_{n-1,j}\}$,
for $1 \le j \le (n-1)$.

One  needs to control what is going on at infinity.  As a matter of
fact, outside a specific rectangle  centered at the origin, the zero set
is the diagonal $\{x=y\}$, see Proposition~\ref{P-ext}.

Hence in this way (like for the square), we get that inside the zero
set, we have the diagonal and $\frac{n-1}{2}$  disjoint simple
closed lines turning around the origin.

\subsection{No other closed curve in the nodal set of
$\Phi_n^{\frac{3\pi}{4}}$}\label{nocc}

It remains to show that there are no other closed curves which do
not cross the diagonal. The ``energy'' considerations of our
previous papers work in the following way.

Assume there is a connected component of the nodal set which does
not meet the lattice $\cL$. Using Proposition~\ref{P-ext}, we see
that this component must be contained in some large coordinate
square centered at $(0,0)$, call it $C$. Since the nodal set cannot
meet the vertical or horizontal lines defined by the zeros of $H_n$,
we would have a nodal domain $\omega$ contained in $C$, hence in one
of the bounded connected components of $\{H_n(x)\, H_n(y) < 0\}$,
and hence also in some infinite rectangle $R$ between two
consecutive zeros of $H_n$. We can compute the energy for $\omega$
by applying Green's formula in $\omega$. We can compute the energy
of the infinite rectangle $R$ by applying Green's formula first in a
finite rectangle, and then taking the limit (using the decaying
exponential factor). We have that the first Dirichlet eigenvalues
$\lambda_1$ satisfy $\lambda_1(\omega) = \lambda_1(R) = 2(n+1)$. On
the other hand, taking some $\omega_1$ such that $\omega \subset
\omega_1 \subset R$, with strict inclusions, we have
$\lambda_1(\omega) > \lambda_1(\omega_1) \ge \lambda_1(R)$, a
contradiction.

A simple alternative argument is the following. We look at
the line $y=\alpha x $ for some $\alpha \neq 1$. The intersection of
the zero set with this line corresponds to the zeroes of the
polynomial $x\mapsto H_n(x) - H_n (\alpha x)$ which has at most $n$
zeroes. But in our previous construction, we get at least $n$
zeroes. So the presence of extra  curves would lead to a
contradiction for some $\alpha$. This  argument  solves the
problem at  infinity as well.

\subsection{Perturbation argument}

\begin{figure}[!ht]
\begin{center}
\includegraphics[width=12cm]{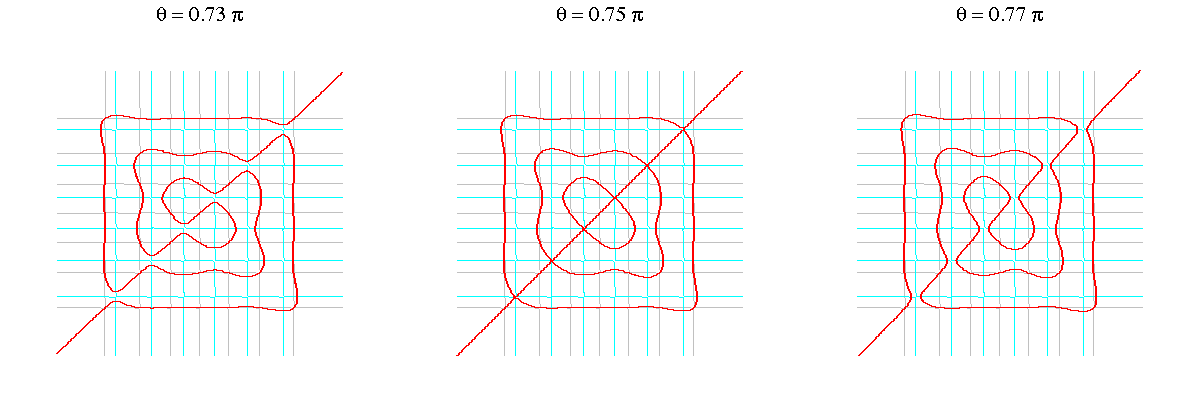}
\caption{The nodal set of $N(\Phi^{\theta}_n)$ near
$\frac{3\pi}{4}$.} \label{FSch-2}
\end{center}
\end{figure}

 Figure~\ref{FSch-2} shows the desingularization of the nodal
set $N(\Phi^{\frac{3\pi}{4}}_n)$, from below and from above. The
picture is the same as in the case of the square (see Figure
\ref{sq-1}), all the critical points disappear at the same time and
in the same manner, \ie all the double crossings open up
horizontally or vertically depending whether $\theta$ is less than
or bigger than $\frac{3\pi}{4}$.

\begin{figure}[!ht]
\begin{center}
\includegraphics[width=8cm]{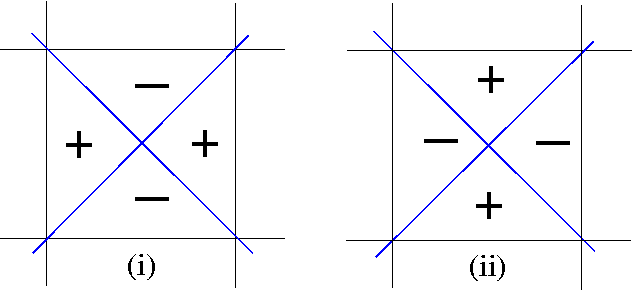}
\caption{Signs near the critical zeros.} \label{FSchr-local-signs-1}
\end{center}
\end{figure}

 As in the case of the square, in order to show that the nodal
set can be desingularized under small perturbation, we look at the
signs of the eigenfunction $\Phi_n^{\frac{3\pi}{4}}$ near the
critical zeros. We use the cases [i] and [ii] which appear in
Figure~\ref{FSchr-local-signs-1} below (see also
\cite[Figure~6.7]{BeHe}).

The sign configuration for $\phi_{n,0}(x,y) - \phi_{0,n}(x,y)$ near
the critical zero $(t_{n-1,i},t_{n-1,i})$ is that of
\begin{equation*}\label{sign}
\left\lbrace
\begin{array}{l}
\text{case [i], ~if~} i \text{~is even},\\
\text{case [ii],~if~} i \text{~is odd}.\\
\end{array}\right.
\end{equation*}

Looking at the intersection of the nodal set with the vertical line
$\{y=t_{n-1,i}\}$, we have that
$$(-1)^i \left( H_n(t) - H_n(t_{n-1,i})\right) \ge 0, \text{~for~}
t \in ]t_{n,i},t_{n,i+1}[ \,.
$$

For positive $ \epsilon$ small, we write

$$(-1)^i \left( H_n(t) - (1+\epsilon) H_n(t_{n-1,i})\right) =
(-1)^i \left( H_n(t) - H_n(t_{n-1,i})\right) + \epsilon (-1)^{i+1}\,
H_n(t_{n-1,i})\,,
$$

so that

$$(-1)^i \left( H_n(t) - (1+\epsilon) H_n(t_{n-1,i})\right) \ge 0,
\text{~for~} t \in ]t_{n,i},t_{n,i+1}[ \,.
$$

\begin{figure}[!ht]
\begin{center}
\includegraphics[width=8cm]{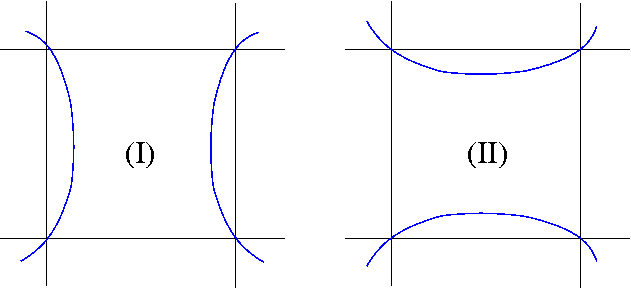}
\caption{Desingularization at a critical zero.}
\label{FSchr-local-signs-2}
\end{center}
\end{figure}

 A similar statement can be written for horizontal line
$\{x=t_{n-1,i}\}$ and $-\epsilon$, with $\epsilon > 0$, small
enough. These inequalities describe how the crossings open up all at
the same time, and in the same manner, vertically (case I) or
horizontally (case II), see Figure~\ref{FSchr-local-signs-2}, as in
the case of the square \cite[Figure~6.8]{BeHe}.

We can then conclude as in the case of the square, using the local
nodal patterns, Section~\ref{SS-lnp}.\medskip

 \textbf{Remark}. Because the local nodal patterns can only
change when $\theta$ passes through one of the values $\theta(i,j)$
defined in \eqref{P10}, the above arguments work for $\theta \in
J\setminus \{\frac{3\pi}{4}\}$, for any interval $J$ containing
$\frac{3\pi}{4}$ and no other value $\theta(i,j)$.

\section{Proof of Theorem \ref{th2}}

\begin{proposition}\label{P-1}
The conclusion of Theorem \ref{th2} holds with
\begin{equation}\label{P26}
\theta_c := \inf \left\{ \theta(i,j) ~|~ 1 \le i,j \le n-1 \right\}
\,.
\end{equation}
\end{proposition}

\textbf{Proof}. The proof consists in the following steps. For
simplicity, we call $N$ the nodal set $N(\Phi^{\theta}_n)$.
\vspace{-3mm}
\begin{itemize}
    \item Step 1.~ By Proposition~\ref{P-ext}, the structure of the nodal
    set $N$ is known  outside a large coordinate rectangle centered at $(0,0)$
    whose sides are defined by the ad hoc numbers in
    Lemmas~\ref{L-ext1} and \ref{L-ext2}. Notice that the sides of
    the rectangle  serve as barriers for the arguments using the local
    nodal patterns as in our paper for the square.
    \item Step 2.~ For $1\le j \le n-1$, the line $\{x=t_{n-1,j}\}$
    intersects the set $N$ at exactly one point $(t_{n-1,j}, y_j)$,
    with $y_j > t_{n,n}$ when $j$ is odd, \resp with $y_j < t_{n,1}$
    when $j$ is even. The proof is given below, and is similar to the proofs of
    Lemmas~\ref{L-ext1} or \ref{L-ext2}.
    \item Step 3.~ Any connected component of $N$
    has at least one point in common with the set $\cL$. This
    follows from the  argument with $y = \alpha x$ or from the energy argument
    (see Subsection \ref{nocc}).
    \item Step 4.~ Follow the nodal set from the point
    $(t_{n,1},t_{n,n})$ to the point $(t_{n,n},t_{n,1})$, using the
    analysis of the local nodal patterns as in the case of the
    square.
\end{itemize}

\textbf{Proof of Step~2}. For $1 \le j \le (n-1)$, define the
function $v_j$ by
$$
v_j(y) := \cos\theta \, H_n(t_{n-1,j}) + \sin\theta \, H_n(y)\,.
$$

The local extrema of $v_j$ are achieved at the points $t_{n-1,i}$,
for $1 \le i \le (n-1)$, and we have
$$
v_j(t_{n-1,i} ) = \cos\theta \, H_n(t_{n-1,j}) + \sin\theta \,
H_n(t_{n-1,i})\,,
$$
which can be rewritten, using \eqref{P10}, as
$$
v_j(t_{n-1,i}) = \frac{H_n(t_{n-1,j})}{\sin\theta(j,i)}\,
\sin\left(\theta(j,i) -\theta \right).
$$

The first term in the right-hand side has the sign of $(-1)^{j+1}$
and the second term is positive provided that $0 < \theta <
\theta_c$. Under this last assumption, we have
\begin{equation}\label{P-32}
(-1)^{j+1}\, v_j(t_{n-1,i}) > 0, ~\forall i, ~1 \le i \le (n-1)\,.
\end{equation}

The assertion follows. \hfill $\square$

\section{Courant's theorem for the $2D$ quantum harmonic
oscillator}

Recall that $\cE_{\ell}$ is the eigenspace of $\hH$
associated with the eigenvalue $\hat{\lambda}(\ell) := 2(\ell + 1)$.
This eigenspace is generated by the eigenfunctions
$\phi_{\ell-j,j}$, for $0 \le j \le \ell$. It has dimension $(\ell +
1)$. The functions in $\cE_{\ell}$ are even (\resp odd) under the
map $a : (x,y) \mapsto (-x,-y)$ when $\ell$ is even (\resp odd).

Since $\dim(\bigoplus_{j=0}^{\ell-1}\cE_j) = \frac{\ell(\ell
+1)}{2}$, Courant's theorem gives the following estimate for the
number $\mu(u)$ of nodal domains of an eigenfunction $u \in
\cE_{\ell}$,
\begin{equation}\label{courant1}
\mu (u) \le \frac{\ell(\ell +1)}{2} + 1 =: \mu_C(\ell)\,.
\end{equation}

Using the symmetry or anti-symmetry with respect to $a$, one can
improve Courant's estimate.

\begin{proposition}\label{P-courant2}
Let $u \in \cE_{\ell}$. Then, the number $\mu(u)$ of nodal domains
of $u$ satisfies the inequalities,
\begin{equation}\label{courant2}
\mu(u) \le \mu_L(\ell) := \left\lbrace
\begin{array}{ll}
2(r^2+1) & \text{~if~} \ell = 2r\,,\\
2r(r+1)+2 & \text{~if~} \ell = 2r+1\,.\\
\end{array}
\right.
\end{equation}
In particular, we have that $\mu_L(\ell) < \mu_C(\ell)$ provided
that $\ell \ge 3$, and $\mu_L(\ell) = \mu_C(\ell)$, when $\ell = 2$.
\end{proposition}

\begin{corollary}\label{C-courant2}
The only Courant sharp eigenvalues of the quantum harmonic
oscillator are the eigenvalues,
\begin{equation}\label{courant-s}
\left\lbrace
\begin{aligned}
& \hat{\lambda}(0) = 2, \text{~with~} \mu_C(0)=1\,,\\
& \hat{\lambda}(1) = 4, \text{~with~} \mu_C(1)=2\,,\\
& \hat{\lambda}(2) = 6, \text{~with~} \mu_C(2)=4\,.\\
\end{aligned}
\right.
\end{equation}
\end{corollary}

\textbf{Proof of the corollary}. The first two assertions are clear.
The last one follows from the fact than the nodal set of an
eigenfunction in $\cE_2$ is a hyperbola, the union of two lines
which intersect, or an ellipse. \hfill $\square$

\textbf{Proof of the proposition}. We use Leydold's argument in
\cite{Ley}, namely the symmetry properties of the eigenfunctions
with respect to $a$, the fact that an odd eigenfunction is always
orthogonal to an even one, and Courant's proof.

\noib Assume that $u \in \cE_{\ell}$ with $\ell = 2r$. We have
$$\dim(\bigoplus_{j=0}^{r-1}\cE_{2j}) = r^2\,.
$$
There are $k_i$ nodal domains of $u$ which are $a$ invariant,
$a(\omega)=\omega$, and $2 k_a$ nodal domains which are not
invariant, $a(\omega) \cap \omega = \emptyset$. Assume that $k_i+k_a
\ge r^2+2$. Define functions $u_j$ such that $u_j = u|_{\omega_j}$,
and $0$ elsewhere, for each invariant domain $\omega_j$, $1 \le j
\le k_i-1$, and $u_j=u|_{\omega_p \cup a(\omega_p)}$, and $0$
elsewhere, for the $k_a$ non-invariant domains. This gives us
$k_a+k_i-1 \ge r^2+1$ independent functions. We can find a linear
combination $v$ of these functions such that $\|v\|_{L^2} = 1$, $v
\perp \bigoplus_{j=0}^{r-1}\cE_{2j}$, and $\mathcal{Q}(v)) = 2(\ell
+1)$, where $\mathcal{Q}$ is the quadratic form associated with
$\hH$. The function $v$ is even by construction so that it is
orthogonal to any odd eigenfunction. It follows that $v \in
\cE_{\ell}$ which leads to a contradiction since $v$ vanishes on an
open set. It follows that $k_i+k_a \le r^2+1$ and hence that $k_i+
2k_a \le 2(r^2+1)$. This proves the first assertion.

\noib Assume that $u \in \cE_{\ell}$ with $\ell = 2r+1$. The proof
is similar. We have that $$\dim(\bigoplus_{j=0}^{r-1}\cE_{2j+1}) =
r(r+1)\,.$$
For an odd eigenfunction $u$, the nodal domains satisfy
$a(\omega) \cap \omega = \emptyset$, so that $\mu(u) = 2k$, and we
can construct $k$ linearly independent functions $u_j=u|_{a(\omega)
\cup \omega}$, and we can proceed as above. \hfill $\square$

 \textbf{Remark}. In the above proof, we used Courant's proof
which is based on energy estimates, using Green's formula for the
eigenfunction $u$. That this can be done in the case of the quantum
harmonic operator follows from the following argument. At infinity,
the nodal set of $u$ is a regular submanifold. It consists of arcs
asymptotic to lines determined by the homogeneous higher order terms
in $\exp(\frac {x^2+y^2}{2})\, u(x,y))$. We can apply Green's formula to
the intersections of the nodal domains of $u$ with balls $B(0,r)$.
When $r$ tends to infinity, the boundary terms involving the ball
tend to zero due to the presence of the exponential factors.

\end{document}